\documentclass[12pt]{amsart}
\newtheorem{thm}{Theorem}
\newtheorem{prop}[thm]{Proposition}

\newtheorem{cor}[thm]{Corollary}

\theoremstyle{remark}
\newtheorem{rem}[thm]{Remark}

\theoremstyle{definition}

\newcommand{\FF}{\mathcal{ F}}
\newcommand{\II}{\mathcal{ I}}

\newcommand{\R}{\mathbb{ R}}

\newcommand{\al}{\alpha}
\newcommand{\bb}{\beta}
\newcommand{\cc}{\gamma}
\newcommand{\ww}{\wedge}

\title{Godbillon-Vey invariants for families of foliations}
\author{D.~Kotschick}
\address{Mathematisches Institut, Ludwig-Maximilians-Universit\"at 
M\"unchen, Theresienstr.~39, 80333 M\"unchen, Germany}
\email{dieter@member.ams.org}

\date{\today; MSC 2000: 57R30, 58A17, 58H15, 37C85.}

\thanks{Support from the {\it VolkswagenStiftung} and the hospitality 
of Harvard University are gratefully acknowledged.}

\begin{document}

\begin{abstract}
The classical Godbillon-Vey invariant is an odd degree 
cohomology class that is a cobordism invariant of a single foliation. 
Here we investigate cohomology classes of even degree that are cobordism 
invariants of (germs of) $1$-parameter families of foliations.
\end{abstract}

\maketitle

In this paper we study an analogue of the classical Godbillon-Vey 
invariant~\cite{GV} that is an invariant not of a single foliation 
but of a family of foliations depending smoothly on a real parameter.
For families of foliations of codimension $q$, this invariant is
a family of cohomology classes in degree $2q+2$. It is constructed and
investigated here using explicit elementary calculations with differential 
forms in the spirit of the original work of Godbillon and Vey. We do 
not use Weil algebras or Gelfand-Fuks cohomology.

Recall that associated to a codimension $q$ foliation $\FF$ with defining 
$q$-form $\alpha$ one has a $1$-form $\beta$ such that 
$d\alpha=\alpha\wedge\beta$. This has the remarkable property that 
$\beta\wedge (d\beta)^{q}$ is closed, and that its cohomology class is 
independent of the particular choices made for $\alpha$ and $\beta$. 
This cohomology class is the Godbillon-Vey invariant of $\FF$. At 
least for codimension one, the 
calculations showing the independence of choices are quite concise, 
but non-trivial, see Section~\ref{s:GV} below. 

In the case of a $1$-parameter family $\FF_{t}$, the defining form 
$\alpha$ depends on $t$, and so does $\beta$. Denoting the time 
derivative by a dot, our invariant is the cohomology class $TGV(\FF_{t})$ 
represented by the form $\dot\beta\wedge\beta\wedge (d\beta)^{q}$. It is 
still easy to see that this form is closed. However, proving the 
independence of choices comes down to a highly non-trivial and quite 
miraculous calculation. This is carried out for $q=1$ in 
Section~\ref{s:4D}, and for arbitrary $q$ in Section~\ref{s:higher}.

Gelfand-Feigin-Fuks~\cite{GFF} introduced certain characteristic classes 
of $1$-parameter families of foliations from the point of view of the 
cohomology of the Lie algebra of formal vector fields. Elaborating on 
this, Fuks~\cite{Fuks} mentioned that $\dot\beta\wedge\beta\wedge d\beta$
represents a characteristic class of $1$-parameter families of 
codimension one foliations, without giving a proof that it is 
well-defined. In a recent paper~\cite{lodder}, Lodder tried to supply 
this proof, but his calculation was flawed, in that he disregarded 
several terms involving the time derivatives, see Remark~\ref{r:Lodder}. 
For $q>1$, the forms $\dot\beta\wedge\beta\wedge (d\beta)^{q}$ have not 
been considered in the literature at all.

Our invariant of families admits a factorisation into a product of the 
Reeb class with a cohomology class $T(\FF_{t})$ of degree $2q+1$ in the 
ideal of the foliation. This is analogous to the factorisation of the 
classical Godbillon-Vey invariant discovered by Duminy. In the case of 
codimension one, the class $T(\FF_{t})$ detects an infinitesimal rigidity 
of the Godbillon-Vey classes of transversely homographic foliations 
under arbitrary variations, which are allowed to go outside the subset 
of transversely homographic foliations. For variations only among 
transversely homographic foliations, rigidity of the Godbillon-Vey 
invariant was proved in~\cite{BG}. 

In higher codimension, but not in codimension one, $T(\FF_{t})$ coincides, 
up to a universal constant, with the time-derivative of the classical 
Godbillon-Vey invariant in a family. There is another difference between 
the case of codimension one and that of higher codimension: in 
codimension one the class $TGV(\FF_{t})$ can be interpreted as an 
asymptotic linking number between certain foliations and vector 
fields~\cite{KV}, and this interpretation breaks down in higher 
codimension.

In the first section we recall the definition of the classical 
Godbillon-Vey invariant of a transversely oriented codimension one 
foliation. These well-known calculations will be used later in the 
paper and also serve to fix our notation. In Section~\ref{s:4D} 
we introduce degree $4$ cohomology classes that are invariants of 
families of codimension one foliations, substantiating the claims 
in~\cite{Fuks,lodder}. In Section~\ref{s:vanish} we 
prove several vanishing theorems; in particular we show that the  
invariants vanish for families generated by moving a fixed 
foliation by a flow. We also show that there are many foliations for 
which the invariant vanishes for all (infinitesimal) deformations. 
In Section~\ref{s:dec} we show that the new invariants factorise 
through the Godbillon operator, and discuss the relationship with the 
time-derivative of the classical Godbillon-Vey invariant. In 
Section~\ref{s:higher} we extend the whole discussion to higher 
codimension. Section~\ref{s:last} contains a few remarks.

We work with smooth foliations and smooth families, with smoothness 
of class $C^{\infty}$. For background on the Godbillon-Vey 
invariant we refer to~\cite{Ghys}.

\section{The classical Godbillon-Vey invariant}\label{s:GV}

First, recall the definition of the Godbillon-Vey invariant~\cite{GV}.
Let $\FF$ be a smooth codimension $1$ foliation on a smooth manifold $M$
(of arbitrary dimension). We shall assume throughout that $\FF$
is transversely oriented, so that the normal bundle is trivialised. This
situation can always be achieved by passing to a 2-fold covering.

We denote by $\II (\FF )$ the differential ideal of forms vanishing on 
$T\FF$. The product of any two forms in $\II(\FF)$ vanishes identically. 
Let $\al\in\II (\FF )$ be a defining 1-form for $\FF$; we assume that it is
positive on positively oriented transversals. By the Frobenius theorem, we
have
\begin{equation}\label{eq:da}
d\al = \al\ww\bb \ ,
\end{equation}
for some $\bb$.

\begin{thm}\label{t:gv}
The form $\bb\ww d\bb$ is closed, and its cohomology class is an invariant
of $\FF$. In particular, it is independent of the choices made for
$\al$ and $\bb$.
\end{thm}
\begin{proof}
Differentiating~\eqref{eq:da} and substituting back from
it, we find $\al\ww d\bb =0$, and so
\begin{equation}\label{eq:db}
d\bb = \al\ww\cc \ ,
\end{equation}
for some $\cc$.

Now $d(\bb\ww d\bb )=d\bb\ww d\bb =0$, as $d\bb\in\II(\FF)$ 
by~\eqref{eq:db}. Thus $\bb\ww d\bb$ is closed and defines a
de Rham cohomology class $GV(\FF )\in H^3(M,\R )$. To see that it
is well-defined, suppose first that we fix $\al$, but make a
different choice for $\bb$. We replace $\bb$ by $\bb +\tau$,
for some $\tau\in\II (\FF )$. Then $d\tau , d\bb\in\II (\FF )$, 
and the product of any two forms in $\II (\FF )$ vanishes.
Thus $\bb\ww d\bb$ is replaced by
$$
(\bb +\tau )\ww d(\bb + \tau )=\bb\ww d\bb + d(\tau\ww\bb ) \ ,
$$
which represents the same cohomology class.
Thus, for fixed $\al$, we can use any $\bb$ satisfying~\eqref{eq:da}
in the definition of $GV(\FF )$. It remains to be shown that $GV(\FF )$ is
independent of the choice of $\al$.

The other forms defining $\FF$ (which are also positive on positively
oriented transversals) are of the form $f\al$, with $f$ a positive smooth
function on $M$. We have
$$
d(f\al )=df\ww\al +f\al\ww\bb = \al\ww (-df+f\bb )=
f\al\ww (\bb -d\text{log}f) \ .
$$
Thus $\bb\ww d\bb$ is replaced by
$$
(\bb -d\text{log}f)\ww d(\bb -d\text{log}f)=
\bb\ww d\bb -d((\text{log}f)d\bb ) \ ,
$$
which represents the same cohomology class.
\end{proof}

The Godbillon-Vey invariant $GV(\FF)$ is an invariant of cobordism
classes of foliations, in the following sense:
\begin{prop}\label{p:gvcobord}
Given two oriented codimension one foliations $\FF'$ and $\FF''$ on
closed oriented $3$-manifolds $M'$ and $M''$, suppose that there is an
oriented cobordism $W$ between $M'$ and $M''$, with an oriented 
codimension one foliation $\FF$ which restricts to $\FF'$ and 
$\FF''$
on $M'$ and $M''$ respectively. Then $GV(\FF')=GV(\FF'')\in\R$.
\end{prop}
This follows directly from Stokes's theorem, compare the proof of
Proposition~\ref{p:dgvcobord} below. 

\begin{rem}\label{r:natural}
The obvious naturality of the above calculations under pullbacks 
implies invariance properties also in the case where $M$ has larger 
dimension, so that the Godbillon-Vey invariant is not in the top 
dimension. Suppose for example that we have a foliation $\FF$ on
$M\times [0,1]$ which restricts to $\FF'$ at one end and to $\FF''$
at the other. Then $GV(\FF')=GV(\FF'')\in H^{3}(M,\R)$.
\end{rem}

\section{A $4$-dimensional Godbillon-Vey invariant}\label{s:4D}

We now consider a smoothly varying $1$-parameter family $\FF_t$ of
transversely oriented codimension one foliations on $M$. In this 
section all forms and functions are functions of the parameter 
$t\in\R$. Derivatives with respect to $t$ are denoted by a dot. For
example, let $\al = \al (t)$ be a $1$-parameter family of $1$-forms
with $\al (t)$ defining $\FF_t$. Then equation~\eqref{eq:da} holds,
with $\bb$ also depending on $t$. Let $\dot\bb$ be its time-derivative.

\begin{thm}\label{t:4D}
For every $t\in\R$ the form $(\dot\bb\ww\bb\ww d\bb )(t)$ is closed and 
its cohomology class is an invariant of the family $\FF_t$. In particular,
it is independent of the choices made for $\al$ and $\bb$.
\end{thm}
\begin{proof}
Equation~\eqref{eq:db} holds as before, and $\cc$ also depends on
$t$. Differentiating~\eqref{eq:db} with respect to $t$ we obtain
\begin{equation}\label{eq:dbt}
d\dot\bb=\dot\al\ww\cc + \al\ww\dot\cc \ .
\end{equation}
Combining this with~\eqref{eq:db}, we obtain 
$d\dot\bb\ww d\bb =0$.

Now 
$d(\dot\bb\ww\bb\ww d\bb )=d\dot\bb\ww\bb\ww d\bb -\dot\bb\ww d(\bb\ww d\bb )$,
with the first summand vanishing by the above argument. The second
summand vanishes because $\bb\ww d\bb$ is closed as shown in the proof
of Theorem~\ref{t:gv}. Thus $\dot\bb\ww\bb\ww d\bb$ is also closed. 
To see that its cohomology class is independent of choices, we adapt the
proof of Theorem~\ref{t:gv}.

Suppose that $A, B\in\II (\FF )$. Then $A\ww B=0$, and taking the 
time differential, we find
\begin{equation}\label{eq:AB}
\dot A\ww B+A\ww\dot B=0 \ .
\end{equation}

If we replace $\bb$ by $\bb +\tau$
for some $\tau\in\II (\FF )$, then $d\tau , d\bb\in\II (\FF )$, 
and the product of any two forms in $\II (\FF )$ vanishes.
From the proof of Theorem~\ref{t:gv} above we know that $\bb\ww d\bb$
is replaced by $\bb\ww d\bb +d(\tau\ww\bb )$.
Thus $\dot\bb\ww\bb\ww d\bb$ is replaced by
$$
(\dot\bb +\dot\tau )\ww (\bb\ww d\bb +d(\tau\ww\bb ))=
\dot\bb\ww\bb\ww d\bb +\dot\tau\ww\bb\ww d\bb +\dot\bb\ww d\tau\ww\bb 
+\dot\tau\ww d\tau\ww\bb \ .
$$
Consider now the exact form
$$
d(-\dot\bb\ww\tau\ww\bb -\frac{1}{2}\dot\tau\ww\tau\ww\bb )=
-d\dot\bb\ww \tau\ww\bb +\dot\bb\ww d\tau\ww\bb
-\frac{1}{2}d\dot\tau\ww\tau\ww\bb
+\frac{1}{2}\dot\tau\ww d\tau\ww\bb \ ,
$$
where we have used again that the product of any two forms in 
$\II (\FF )$ vanishes. The second summand on the right-hand side
matches the third summand in the calculation above. 
Using~\eqref{eq:AB}, the last two summands are equal, and their sum 
equals $\dot\tau\ww d\tau\ww\bb$. Using~\eqref{eq:AB} again, we
have
$$
\dot\tau\ww d\bb +\tau\ww d\dot\bb=0 \ ,
$$
and so $\dot\tau\ww\bb\ww d\bb = \bb\ww\tau\ww d\dot\bb =
-d\dot\bb\ww\tau\ww\bb$.

Thus, for fixed $\al$, we can use any $\bb$ satisfying~\eqref{eq:da},
and the cohomology class of $(\dot\bb\ww\bb\ww d\bb )(t)$ will be 
independent of this choice. It remains to be seen that it is
also independent of the choice of $\al$.

The other forms defining $\FF_t$ (which are also positive on positively
oriented transversals) are of the form $f\al$, with $f$ a positive smooth
function on $M$, again depending on $t$. From the proof of Theorem~\ref{t:gv}
we know that replacing $\al$ by $f\al$ results in replacing
$\bb\ww d\bb$ by $\bb\ww d\bb -d((\text{log}f)d\bb )$. Now $\dot\bb$
is replaced by $\dot\bb -d(\frac{d}{dt}\text{log}f)=\dot\bb -d(\frac{1}{f}\dot f)$.

Putting these calculations together, $\dot\bb\ww\bb\ww d\bb$ is replaced
by
\begin{alignat*}{1}
&(\dot\bb-d(\frac{1}{f}\dot f))\ww (\bb\ww d\bb -d((\text{log}f)d\bb ))=\\
&\dot\bb\ww\bb\ww d\bb -d(\frac{1}{f}\dot f)\ww\bb\ww d\bb -\dot\bb\ww
d((\text{log}f)d\bb )+d(\frac{1}{f}\dot f)\ww d((\text{log}f)d\bb )=\\
&\dot\bb\ww\bb\ww d\bb -d(\frac{1}{f}\dot f\ww\bb\ww d\bb) +d(\dot\bb\ww
(\text{log}f)d\bb )+d(\frac{1}{f}\dot f\ww d((\text{log}f)d\bb )) \ ,
\end{alignat*}
which represents the same cohomology class. Here we have used that 
$\bb\ww d\bb$ is closed in order to rewrite the second summand,
and~\eqref{eq:dbt} to rewrite the third one.
\end{proof}

\begin{rem}\label{r:Lodder}
A version of Theorem~\ref{t:4D} appears in~\cite{lodder}. The proof
of well-definedness of the cohomology class of $(\dot\bb\ww\bb\ww d\bb )(t)$ 
is incomplete in~\cite{lodder}, because when replacing $\al$ by $f\al$
and $\bb$ by $\bb +\tau = \bb + g\al$, the author ignores the terms involving 
the time derivatives of $f$ and $g$.
\end{rem}

Because of Theorem~\ref{t:4D}, the family of cohomology classes 
$TGV(\FF_t)=[(\dot\bb\ww\bb\ww d\bb )(t)]\in H^4(M,\R )$ is an invariant 
of the family $\FF_t$. This can be specialised to a single cohomology
class either by integrating over $t$, say
\begin{equation}
IGV(\FF_t )=\int_0^1TGV(\FF_t)dt \ ,
\end{equation} 
or by evaluating at a specific value of $t$, say
\begin{equation}
DGV(\FF_t )=TGV(\FF_t)(0) \ .
\end{equation} 
It is clear that $DGV(\FF_t )$ is an invariant of the germ of $\FF_t$ at 
$t=0$, so it does not depend on the whole family of foliations $\FF_t$.
In fact, any $1$-form $\omega$ with 
$$
\omega\ww d\al + \al\ww d\omega =0
$$
can be considered as an infinitesimal variation of the foliation defined
by $\al$, and in the case that this infinitesimal variation integrates to
an actual variation with $\dot\al = \omega$, we have 
$$
\dot\beta\ww\bb\ww d\bb = \dot\bb\ww\bb\ww\al\ww\cc = d\omega\ww\bb\ww\cc \ ,
$$
so that the right-hand side can be used to define $DGV$ for any infinitesimal
variation of a codimension one foliation.
 
The family of $4$-dimensional cohomology classes $TGV(\FF_t )$ is a cobordism 
invariant of families in the following sense, cf.~Proposition~\ref{p:gvcobord}:
\begin{prop}\label{p:dgvcobord}
Given two smooth families of oriented codimension one foliations 
$\FF'_t$ and $\FF''_t$ on closed oriented 4-manifolds $M'$ and $M''$, 
suppose that there is an oriented cobordism $W$ between $M'$ and $M''$, 
with a smooth family of oriented codimension one foliations $\FF_t$ 
which restricts to $\FF'_t$ and $\FF''_t$ on $M'$ and $M''$ respectively. 
Then $TGV(\FF'_t)=TGV(\FF''_t)\in\R$.
\end{prop}
\begin{proof}
Let $\al$ be a time-dependent $1$-form defining $\FF_t$ 
on $W$, with $d\al = \al\ww\bb$. Then 
\begin{alignat*}{1}
\langle TGV(\FF''_t),[M'']\rangle &- \langle TGV(\FF'_t),[M']\rangle \\
&= \int_{M''}\dot\bb\ww\bb\ww d\bb (t) - \int_{M'}\dot\bb\ww\bb\ww d\bb (t)\\ 
&= \int_W d(\dot\bb\ww\bb\ww d\bb (t)) =\int_W 0 = 0 \ ,
\end{alignat*}
by Stokes's theorem, because $\dot\bb\ww\bb\ww d\bb (t)$ is closed.
\end{proof}
{\it A fortiori}, $IGV$ and $DGV$ are also cobordism invariants.

{\it Mutatis mutandis}, Remark~\ref{r:natural} applies to $TGV(\FF_{t})$ 
as well.

\section{Vanishing theorems}\label{s:vanish}

It is clear from the definition that $TGV(\FF_t )$
vanishes if $\FF_t$ can be defined 
by a form $\al$ with $d\al =\al\ww\bb$ and $\bb$ closed.
In particular this holds if $\al$ can be chosen to
be closed itself.

This observation can be generalised quite a bit. Recall
equations~\eqref{eq:da} and~\eqref{eq:db}. Applying $d$ to
the latter we find $0=\al\ww (\bb\ww\cc -d\cc )$. Therefore,
there exists a $1$-form $\delta $, such that
\begin{equation}\label{eq:dc}
d\cc =\bb\ww\cc +\al\ww\delta \ .
\end{equation}
Repeating the procedure, we also find that
\begin{equation}\label{eq:dd}
d\delta =2\bb\ww\delta +\al\ww\epsilon \ ,
\end{equation}
for some $\epsilon$.

With this notation, we have the following vanishing theorem:
\begin{thm}\label{t:vanish}
If $\FF_t$ can be defined by a $1$-form $\al$ for which
$\bb$ or $\cc$ or $\delta$ vanishes, then $TGV(\FF_t )=0$.

If $\FF_t$ can be defined by a $1$-form $\al$ such that $\al$
or $\bb$ or $\cc$ or $\delta$ is closed, then $TGV(\FF_t )=0$.
\end{thm}
\begin{proof}
We use the above equations and the
time-differential of~\eqref{eq:da}. First notice
\begin{alignat*}{1}
\dot\bb\ww\bb\ww &d\bb=\dot\bb\ww\bb\ww\al\ww\cc =d\dot\al\ww\bb\ww\cc \\
&= d(\dot\al\ww\bb\ww\cc )-\dot\al\ww\bb\ww d\cc = d(\dot\al\ww\bb\ww\cc )
-\dot\al\ww\bb\ww\al\ww\delta \ .
\end{alignat*}
Thus $TGV(\FF_t )=0$ as soon as $\bb$ or $\cc$ or $\delta$ 
vanishes identically. The vanishing of $\bb$, respectively $\cc$,
is equivalent to $\al$, respectively $\bb$, being closed.
The above calculation also shows that $TGV(\FF_t )$ vanishes
if $\cc$ is closed. 

It remains to see what happens when $\delta$ is closed. We 
continue with the last term in the above calculation and 
use~\eqref{eq:dd}:
$$
\dot\al\ww\bb\ww\al\ww\delta = 
-\dot\al\ww\al\ww\frac{1}{2}(d\delta -\al\ww\epsilon )
= -\frac{1}{2}\dot\al\ww\al\ww d\delta \ .
$$
Thus $TGV(F_t )$ is represented by an exact form if $\FF_t$
can be defined by a form $\al$ for which $\delta$ can be taken 
to be closed.
\end{proof}
Evaluating at $t=0$ we obtain from the same arguments:
\begin{thm}\label{t:germvanish}
If $\FF_0$ can be defined by a $1$-form $\al$ for which
$\bb$ or $\cc$ or $\delta$ vanishes, then $DGV(\FF_t )=0$
for all germs $\FF_t$ at $\FF_0$.

If $\FF_0$ can be defined by a $1$-form $\al$ such that $\al$
or $\bb$ or $\cc$ or $\delta$ is closed, then $DGV(\FF_t )=0$
for all germs $\FF_t$ at $\FF_0$.
\end{thm}

Foliations with a defining form $\al$ for which $\delta$
can be taken to be zero have special transverse structures.
These are the {\sl transversely homographic} structures 
in the sense of~\cite{G}, page 174. If $\omega$ is a
$1$-form considered as an infinitesimal variation of $\FF_0$,
a calculation like the one in the proof of Theorem~\ref{t:vanish}
shows that the derivative of $GV$ at $\FF_0$
in the direction of $\omega$ is the cohomology class of
$-2\omega\ww\al\ww\delta$. Thus transversely homographic
foliations are the critical points of $GV$, and 
for these Theorem~\ref{t:germvanish} gives:
\begin{cor}\label{c:homogr}
Let $\FF_0$ be a codimension one foliation with a transversely
homographic structure. Then for all germs $\FF_t$ based at $\FF_0$
we have $DGV(\FF_t )=0$.
\end{cor}

Another byproduct of the proof of Theorem~\ref{t:vanish} is the 
observation that variations in the direction of a closed
$1$-form, i.~e.~with $d\dot\al=0$, always give $TGV(F_t )=0$.
%In this case $\dot\al$ is in the span of $\al$ and $\bb$.
A similar argument shows that $TGV(\FF_t )$ vanishes for a family
of diffeomorphic foliations:
\begin{thm}\label{t:diff}
Let $\FF$ be a smooth codimension one foliation on $M$, and $\Phi_t$ 
a smooth 1-parameter family of diffeomorphisms of $M$. If $\FF_t$ is 
the family generated by pulling back $\FF$ via $\Phi_t$, 
then $TGV(\FF_t )=0$.
\end{thm}
\begin{proof}
By definition, we can choose $\al$ so that $\dot\al =L_X\al$,
with $X$ a (time-dependent) vector field. Now we calculate:
\begin{alignat*}{1}
\dot\bb\ww\bb\ww d\bb &=\al\ww\dot\bb\ww\bb\ww\cc
=(d\dot\al)\ww\bb\ww\cc - \dot\al\ww\bb\ww\bb\ww\cc\\
&=(dL_X\al )\ww\bb\ww\cc = (di_Xd\al )\ww\bb\ww\cc\\
&=-d(\bb (X))\ww\al\ww\bb\ww\cc = d(\bb (X))\ww\bb\ww d\bb\\
&=d(\bb (X)\bb\ww d\bb ) \ .
\end{alignat*}
Here the step from the second line to the third is achieved
by the following:
\begin{alignat*}{1}
di_Xd\al &=di_X(\al\ww\bb )=d(\al (X)\bb -\bb (X)\al )\\
&=d(\al (X))\ww\bb +\al (X)d\bb -d(\bb (X))\ww\al -\bb (X)d\al\\
&=d(\al (X))\ww\bb +\al (X)\al\ww\cc -d(\bb (X))\ww\al -\bb (X)\al\ww\bb \ .
\end{alignat*}

Thus $TGV(\FF_t )$ is represented by an exact form.
\end{proof}

\section{Factoring through the Reeb class}\label{s:dec}

There is a factorisation of the classical Godbillon-Vey invariant
into two different invariants due to Duminy, see~\cite{Ghys}. We 
shall explain in this section that there is an analogous factorisation
of the invariant $TGV(\FF_{t})$.

Recall that $\II(\FF)$ is a graded differential ideal in the algebra 
$\Omega^{*}(M)$ of smooth differential forms on $M$. We shall denote
its cohomology by $H^{*}(\II(\FF))$.

The quotient $\Omega^{*}(\FF)$ of $\Omega^{*}(M)$ by $\II(\FF)$
can be thought of as differential forms on $M$ defined
only along the leaves of $\FF$, with the induced differential 
$d_{\FF}$ being differentiation along the leaves. We shall 
write $H^{*}(\FF)$ for the cohomology of the complex 
$(\Omega^{*}(\FF),d_{\FF})$.

Note that $H^{*}(\II(\FF))$ is a module over the de Rham 
cohomology of $M$. As the product of any two forms in 
$\II(\FF)$ vanishes, we obtain a bilinear map
\begin{equation}\label{eq:cup}
H^{i}(\II(\FF))\times H^{j}(\FF)\longrightarrow H^{i+j}(\II(\FF))
\end{equation}
for all $i$ and $j$.

Suppose now that $\FF$ is defined by $\al$ with $d\al=\al\ww\bb$.
Then $d\bb\in\II(\FF)$; equivalently the projection of $\bb$ to
$\Omega^{1}(\FF)$ is $d_{\FF}$-closed. Thus $\bb$ defines a 
cohomology class $[\bb]\in H^{1}(\FF)$. We saw in the proof of
Theorem~\ref{t:gv} that making a different choice for $\bb$, 
respectively for $\al$, will change $\bb$ by an element of $\II(\FF)$,
respectively by an exact form. Thus the {\it Reeb class} $[\bb]\in H^{1}(\FF)$
is independent of these choices and is an invariant of the 
foliation $\FF$.

Similarly, $d\bb$ is in the ideal $\II(\FF)$ and is obviously
closed. The proof of Theorem~\ref{t:gv} shows that making 
different choices changes $d\bb$ at most by the addition of $d\tau$
with $\tau\in\II(\FF)$. Thus the cohomology class $[d\bb]\in 
H^{2}(\II(\FF))$ is well-defined. Taking the product of this 
class with the Reeb class according to~\eqref{eq:cup} yields 
a class in $H^{3}(\II(\FF))$ whose image in the de Rham cohomology
is the classical Godbillon-Vey invariant $GV(\FF)$.

Here is the analogous statement for $TGV(\FF_{t})$.

\begin{thm}\label{t:dv}
Let $\FF_{t}$ be a smooth family of smooth codimension one 
foliations on $M$. Then for every $t$, the $3$-form 
$(\dot\bb\ww d\bb)(t)$ is in $\II(\FF_{t})$. It is closed, and its 
cohomology class $T(\FF_{t})\in H^{3}(\II(\FF_{t}))$ is independent 
of choices. Its image in the de Rham cohomology is 
$\frac{1}{2}\frac{d}{dt}GV(\FF_{t})$.

The product of $T(\FF_{t})$ with the Reeb class according 
to~\eqref{eq:cup} is a class in $H^{4}(\II(\FF))$ whose image
in the de Rham cohomology equals $TGV(\FF_{t})$ up to sign.
\end{thm}
\begin{proof}
As $d\bb$ is in the ideal, so is $\dot\bb\ww d\bb$. We have 
$d(\dot\bb\ww d\bb)=d\dot\bb\ww d\bb$. We noted in the proof of
Theorem~\ref{t:4D} that this vanishes by the 
combination of~\eqref{eq:db} and~\eqref{eq:dbt}.

That $[\dot\bb\ww d\bb]\in H^{3}(\II(\FF_{t}))$ is independent of 
choices also follows from the proof of Theorem~\ref{t:4D}. In more
detail, suppose that we replace $\bb$ by $\bb+\tau$ with
$\tau\in\II(\FF_{t})$. Then $\dot\bb\ww d\bb$ is replaced by
$\dot\bb\ww d\bb -d(\dot\bb\ww\tau+\frac{1}{2}\dot\tau\ww\tau)$,
where $\dot\bb\ww\tau$ and $\dot\tau\ww\tau$ are in $\II(\FF_{t})$
because $\tau$ is. If we replace $\al$ by $f\al$, then
$\dot\bb\ww d\bb$ is replaced by $\dot\bb\ww d\bb -d(\frac{1}{f}\dot f 
d\bb)$. Thus $T(\FF_{t})$ is well-defined in $H^{3}(\II(\FF_{t}))$.

To calculate the image of $T(\FF_{t})$ in the de Rham cohomology,
consider the time-derivative of $\bb\ww d\bb$:
\begin{equation}\label{eq:dGV}
\frac{d}{dt}(\bb\ww d\bb)=\dot\bb\ww d\bb+\bb\ww d\dot\bb=
2\dot\bb\ww d\bb+d(\dot\bb\ww\bb) \ .
\end{equation}
The left-hand side represents $\frac{d}{dt}GV(\FF_{t})$
in the de Rham cohomology, and the right-hand side represents
twice the image of $T(\FF_{t})$.

Given the above, the last claim in the Theorem is obvious.
\end{proof}
Note that there is no reason to expect the time-derivative of 
$GV(\FF_{t})$ to define a class in $H^{3}(\II(\FF_{t}))$. The 
formula~\eqref{eq:dGV} shows that, up to a factor of $2$, $T(\FF_{t})$ 
lifts $\frac{d}{dt}GV(\FF_{t})$ to $H^{3}(\II(\FF_{t}))$. But 
the difference term $d(\dot\bb\ww\bb)$ is not usually in the 
ideal $\II(\FF_{t})$. Thus $\dot\bb\ww\bb$ does not define a 
cohomology class in $H^{2}(\FF_{t})$, and the formula~\eqref{eq:dGV} 
does not have cohomological meaning, beyond expressing the 
relationship between the image of $T(\FF_{t})$ in the de Rham 
cohomology and the time derivative of the classical Godbillon-Vey 
invariant. At this point there is a crucial difference between the 
codimension $1$ case and that of higher codimension; see 
Theorem~\ref{t:dvhigher} below.

It is well-known that there are smooth families of foliations
for which the classical Godbillon-Vey invariant is not constant.
Every such family gives examples for the non-vanishing of $T(\FF_{t})$.

This contrasts with the following vanishing result which sharpens 
Corollary~\ref{c:homogr}:
\begin{thm}\label{t:homo}
Let $\FF_0$ be a codimension one foliation with a transversely
homographic structure. Then for all germs $\FF_t$ based at $\FF_0$
we have $T(\FF_t)(0)=0\in H^{3}(\II(\FF_{0}))$.
\end{thm}
\begin{proof}
We calculate with the defining form of $T(\FF_t)$ using~\eqref{eq:db},
the time differential of~\eqref{eq:da}, and then~\eqref{eq:dc} 
and~\eqref{eq:dbt}:
\begin{alignat*}{1}
\dot\bb\ww d\bb 
&=\dot\bb\ww\al\ww\cc=-(d\dot\al)\ww\cc+\dot\al\ww\bb\ww\cc\\
&=-d(\dot\al\ww\cc)-\dot\al\ww d\cc+\dot\al\ww(d\cc-\al\ww\delta)\\
&=d(\al\ww\dot\cc)+\al\ww\dot\al\ww\delta \ .
\end{alignat*}
This shows that $T(\FF_t)\in H^{3}(\II(\FF_{t}))$ is represented by the 
form $\al\ww\dot\al\ww\delta$. By definition, this vanishes for a transversely 
homographic foliation, no matter what $\dot\al$ is.
\end{proof}
As a consequence of the above calculation, note that if a $1$-form
$\omega$ is considered as an infinitesimal variation of a foliation 
$\FF$, we can define the class $T\in H^{3}(\II(\FF))$ for this 
infinitesimal variation as the cohomology class of 
$\al\ww\omega\ww\delta$, even if $\omega$ does not integrate to an 
actual $1$-parameter variation of $\FF$.

It is interesting to examine the other vanishing theorems in 
Section~\ref{s:vanish} in the light of the decomposition of $TGV$. 
If $\bb=0$, then both $T(\FF_{t})$ and the Reeb class vanish. If 
$\cc=0$, then the Reeb class lifts to $H^{1}(M,\R)$, and $T(\FF_{t})=0$.
If $\delta=0$, we can say nothing about the Reeb class, but $T(\FF_{t})$
vanishes by the Theorem above. In Theorem~\ref{t:vanish} we saw that 
$TGV$ vanishes also when $\cc$ or $\delta$ is closed. These results do 
not seem to come from either the vanishing of the Reeb class, or the 
vanishing of $T(\FF_{t})$, but instead rely on the interplay between the 
two via~\eqref{eq:cup}.

Theorem~\ref{t:diff} does in fact come from the vanishing of 
$T(\FF_{t})$ in families generated from a fixed foliation by a flow. 
As a flow acts trivially on the de Rham cohomology, Theorem~\ref{t:dv} 
shows that the image of $T(\FF_{t})$ in the de Rham cohomology is 
trivial. It does not, however, show that $T(\FF_{t})$ vanishes in 
$H^{3}(\II(\FF))$, which itself varies with $t$. Nevertheless, an 
easy adaptation of the proof of Theorem~\ref{t:diff} shows that 
$T(\FF_{t})$ vanishes in the cohomology of the ideal. We shall give 
this argument for foliations of arbitrary codimension in the next 
section.

Dualising the decomposition in Theorem~\ref{t:dv}, we obtain:
\begin{thm}\label{t:dual}
Let $\FF_{t}$ be a smooth family of smooth codimension one 
foliations on a closed oriented manifold $M$. Then for every $t$, 
the class $TGV(\FF_{t})\in H^{4}(M)$, thought of as a linear 
functional on $H^{n-4}(M)$, decomposes up to sign into the 
composition of the following two maps:
\begin{enumerate}
\item the map given by the product with $T(\FF_{t})$:
\begin{alignat*}{1}
    H^{n-4}(M) &\longrightarrow H^{n-1}(\II(\FF_{t}))\\
    x &\longmapsto x\cup T(\FF_{t}) \ ,
\end{alignat*}
\item the Godbillon operator given by the Reeb class according 
to~\eqref{eq:cup}:
\begin{alignat*}{1}
    H^{n-1}(\II(\FF_{t})) &\longrightarrow H^{n}(M,\R)=\R\\
    y &\longmapsto y\cup [\bb] \ .
\end{alignat*}
\end{enumerate}
\end{thm}
The definition of the Godbillon operator above uses~\eqref{eq:cup} 
together with the surjection $H^{n}(\II(\FF))\rightarrow H^{n}(M,\R)$ 
in the top dimension.

This decomposition is very useful because the study of the classical
Godbillon-Vey invariant has led to many vanishing theorems for the 
Godbillon operator in situations where the Reeb class need not 
vanish, cf.~\cite{Ghys}. These vanishing theorems arise from the 
localisation of the Godbillon operator on saturated sets 
discovered by Duminy. For example, it is known that the Godbillon 
operator vanishes for foliations almost without holonomy and for 
those without resilient leaves. Thus, we obtain another 
vanishing theorem:
\begin{thm}\label{t:loc}
Let $\FF_{0}$ be a foliation whose Godbillon operator vanishes,
for example a foliation almost without holonomy, or without 
resilient leaves. Then $DGV(\FF_{t})$ vanishes for all germs $\FF_{t}$ 
at $\FF_{0}$.
\end{thm}

This leaves open the possibility that there may be a foliation
$\FF_{0}$ with $GV(\FF_{0})=0$, but which admits an infinitesimal
deformation with $DGV(\FF_{t})$ non-zero. The latter condition only 
implies the non-triviality of the Godbillon operator, which does not 
contradict the vanishing of the classical Godbillon-Vey invariant.

\section{Higher codimension}\label{s:higher}

There is a classical Godbillon-Vey invariant for foliations of 
higher codimension $q$ with oriented normal bundles. Such foliations 
are defined by locally decomposable $q$-forms $\al = \al_1\ww\ldots\ww\al_q$
of maximal rank. Again we have 
\begin{equation}\label{eq:dal}
d\al =\al\ww\bb \ ,
\end{equation}
for some $1$-form $\bb$. 

\begin{thm}\cite{GV}\label{t:gvhigher}
The form $\bb\ww (d\bb )^q$ is closed and its cohomology class 
$GV(\FF)\in H^{2q+1}(M,\R)$ is an invariant of $\FF$. In particular, 
it is independent of the choices made for $\al$ and $\bb$.
\end{thm}
Like Theorem~\ref{t:gv}, this is proved by a direct calculation. Instead 
of reproducing this lengthy calculation, we proceed directly to the 
decomposition into two invariants and prove that they are well-defined 
and that one recovers the Godbillon-Vey invariant by composition. We then 
generalise the argument to the case of families.

As in the codimension one case, denote by ${\mathcal I}(\FF )$ the 
graded differential ideal of forms on $M$ which vanish when evaluated 
on tuples of vectors all tangent to $\FF$. These are locally linear 
combinations of decomposable forms for which each summand contains one 
of the $\al_{i}$ as a factor. Note that all $(q+1)$-fold products of 
elements of ${\mathcal I}(\FF )$ vanish. We denote the cohomology 
of $\II(\FF)$ by $H^{*}(\II(\FF))$ and use $H^{*}(\FF)$ for the 
cohomology of the quotient of $\Omega^{*}(M)$ by $\II(\FF)$.

Differentiating~\eqref{eq:dal} and substituting back from
it, we find $\al\ww d\bb =0$, and so
\begin{equation}\label{eq:dbe}
d\bb = \sum_{i=1}^q\al_i\ww\cc_i \in {\mathcal I}(\FF ) \ ,
\end{equation}
for some $\cc_i$. Thus $\bb$ defines a cohomology class in 
$H^{1}(\FF)$, called the {\it Reeb class} of $\FF$. As before, this is 
well-defined independent of choices.

The form $d\bb$ is in $\II(\FF )$ by~\eqref{eq:dbe}. 
It is obviously closed and so defines a cohomology class in 
$H^{2}(\II(\FF))$. Making different choices changes $d\bb$ at most by 
the addition of $d\tau$ for some $\tau\in\II(\FF))$. Thus the 
class $[d\bb]\in H^{2}(\II(\FF))$ is well-defined.

In this case, we cannot compose the Reeb class with the class $[d\bb]$ 
to obtain a well-defined cohomology class of degree three. Though 
$H^{*}(\II(\FF))$ is a module over the de Rham cohomology of $M$, if we 
choose a form on all of $M$ representing a cohomology class in 
$H^{*}(\FF)$, then the wedge product with a closed form in the ideal is 
not necessarily closed, so we do not get a pairing of the 
form~\eqref{eq:cup} in all cases. However, if the form in the ideal 
is itself a $q$-fold product of forms in the ideal, then the resulting 
wedge product is closed because $(q+1)$-fold wedge products of forms 
in the ideal vanish. Thus, the product of $[d\bb]^{q}\in 
H^{2q}(\II(\FF))$ with the Reeb class is a well-defined cohomology 
class in $H^{2q+1}(\II(\FF))$, whose image in the de Rham cohomology 
is the Godbillon-Vey class of Theorem~\ref{t:gvhigher}. This of 
course proves Theorem~\ref{t:gvhigher}.

Now we make the same extension to families as in the 
codimension one case. Suppose $\al$ and $\bb$ above depend 
smoothly on a real parameter $t$, and denote the differential
with respect to $t$ by a dot.

\begin{thm}\label{t:dvhigher}
Let $\FF_{t}$ be a smooth family of smooth codimension $q$ foliations 
on $M$. Then for every $t$, the $(2q+1)$-form $(\dot\bb\ww (d\bb)^{q})(t)$ 
is in $\II(\FF_{t})$. It is closed, and its cohomology class
$T(\FF_{t})\in H^{2q+1}(\II(\FF_{t}))$ is a well-defined invariant of 
the family $\FF_{t}$.

If $q\geq 2$, then $T(\FF_{t})=\frac{1}{q+1}\frac{d}{dt}GV(\FF_{t})\in 
H^{2q+1}(\II(\FF_{t}))$.

The product of $T(\FF_{t})$ with the Reeb class is a well-defined 
cohomology class $TGV(\FF_{t})\in H^{2q+2}(\II(\FF_{t}))$ represented 
by $(\dot\bb\ww\bb\ww (d\bb)^{q})(t)$, up to sign.
\end{thm}
\begin{proof}
The form $\dot\bb\ww (d\bb)^{q}$ is in the ideal $\II(\FF_{t})$ 
because $d\bb$ is. 

Differentiating~\eqref{eq:dbe} with respect to $t$, we obtain
\begin{equation}\label{eq:dbet}
d\dot\bb=\sum_{i=1}^q(\dot\al_i\ww\cc_i + \al_i\ww\dot\cc_i) \ .
\end{equation}
We have $d(\dot\bb\ww (d\bb)^q)=d\dot\bb\ww (d\bb )^q$, which vanishes 
by the combination of~\eqref{eq:dbe} and~\eqref{eq:dbet}.

Now suppose that we replace $\bb$ by $\bb+\tau$, with 
$\tau\in\II(\FF_{t})$. Then we have
%\begin{equation}
\begin{alignat*}{1}
&(\dot\bb+\dot\tau)\ww(d\bb+d\tau)^{q}- \dot\bb\ww(d\bb)^{q} \ = \ \\
&\dot\tau\ww  (d\bb)^{q}
+ \sum_{i=0}^{q-1}\binom{q}{i}(\dot\bb+\dot\tau)\ww(d\bb)^{i}\ww(d\tau)^{q-i}
\ = \ \dot\tau\ww(d\tau)^{q}\\
&+ \sum_{i=0}^{q-1}
\left( \binom{q}{i}\dot\bb\ww(d\bb)^i\ww(d\tau)^{q-i}+\binom{q}{i+1}
\dot\tau\ww(d\bb)^{i+1}\ww(d\tau)^{q-i-1} \right) \ .
\end{alignat*}
%\end{equation}
We have to prove that this is the exterior differential of a form in 
$\II(\FF_{t})$.  

As $(q+1)$-fold products of elements in $\II(\FF_{t})$ vanish, 
equation~\eqref{eq:AB} now takes the form
\begin{equation}\label{eq:ABq}
\dot A_0\ww A_1\ww\dots\ww A_q+ \cdots + A_0\ww \dots\ww A_{q-1}\ww\dot A_q= 0
\end{equation}
for all $A_0,\dots,A_q\in \II(\FF)$.

Applying this to $\tau\ww(d\tau)^{q}=0$, we find
$$
\dot\tau\ww(d\tau)^{q}+q\tau\ww d\dot\tau\ww(d\tau)^{q-1}=0 \ .
$$
This implies 
$$
d(\tau\ww\dot\tau\ww(d\tau)^{q-1})=\dot\tau\ww(d\tau)^{q}
-\tau\ww d\dot\tau\ww(d\tau)^{q-1}=
\frac{q+1}{q}\dot\tau\ww(d\tau)^{q} \ .
$$

Similarly, we find 
$$
q\dot\bb\ww(d\bb)^{q-1}\ww d\tau+
\dot\tau\ww(d\bb)^{q}=qd(\tau\ww\dot\bb\ww(d\bb)^{q-1}) \ .
$$

It remains to discuss the terms of the form 
$$
\binom{q}{i}\dot\bb\ww(d\bb)^i\ww(d\tau)^{q-i}+\binom{q}{i+1}
\dot\tau\ww(d\bb)^{i+1}\ww(d\tau)^{q-i-1} \ 
$$
for $i\leq q-2$.
Using~\eqref{eq:ABq} in the form
\begin{alignat*}{1}
\dot\tau\ww(d\bb)^{i+1}\ww(d\tau)^{q-i-1}+ &(i+1)\tau\ww d\dot\bb\ww
(d\bb)^{i}\ww(d\tau)^{q-i-1}\\
+ &(q-i-1)\tau\ww d\dot\tau\ww(d\bb)^{i+1}
\ww(d\tau)^{q-i-2}=0
\end{alignat*}
one easily checks
\begin{alignat*}{1}
&\dot\bb\ww(d\bb)^{i}\ww(d\tau)^{q-i} +\frac{q-i}{i+1}\dot\tau\ww(d\bb)^{i+1}
\ww(d\tau)^{q-i-1} \ = \ \\
&\frac{q-i-1}{i+1}d(\tau\ww\dot\tau\ww(d\bb)^{i+1}\ww
(d\tau)^{q-i-2})-d(\dot\bb\ww\tau\ww(d\bb)^{i}\ww(d\tau)^{q-i-1}) \ ,
\end{alignat*}
where the right-hand side is obviously in $d(\II(\FF_{t}))$.
Because $\binom{q}{i}\frac{q-i}{i+1}= \binom{q}{i+1}$, this finally
proves the independence of the choice of $\bb$.

If we replace $\al$ by $f\al$, then $\bb$ can be replaced by 
$\bb-d\log f$. Thus $\dot\bb\ww(d\bb)^{q}$ is replaced by
$$
\dot\bb\ww(d\bb)^{q}-d(\frac{1}{f}\dot f(d\bb)^{q}) \ ,
$$
which represents the same cohomology class in $H^{2q+1}(\II(\FF_{t}))$.

This completes the proof that $(\dot\bb\ww(d\bb)^{q})(t)$ defines a 
cohomology class $T(\FF_{t})\in H^{2q+1}(\II(\FF_{t}))$ that is independent 
of choices.

Now we calculate the time-derivative of $GV(\FF_{t})$:
\begin{alignat*}{1}
\frac{d}{dt} (\bb\ww(d\bb)^{q}) \ &= \ 
\dot\bb\ww(d\bb)^{q}+q\bb\ww d\dot\bb\ww(d\bb)^{q-1}\\
&= \ (q+1)\dot\bb\ww(d\bb)^{q}+qd(\dot\bb\ww\bb\ww(d\bb)^{q-1}) \ .
\end{alignat*}
If $q\geq 2$, the right-hand side differs from the defining form of 
$(q+1)T(\FF_{t})$ by a form in $d(\II(\FF_{t}))$.

By definition, the defining form of $T(\FF_{t})$ is a $q$-fold product of
elements in the ideal of $\FF_{t}$. Thus, its product with the Reeb 
class is well-defined.
\end{proof}

As in the case of codimension one we can define $IGV(\FF_t )$ by 
integrating $TGV(\FF_{t})$ over $t$ and $DGV(\FF_t )$ by evaluation at $t=0$.

We now prove the generalisation of Theorem~\ref{t:diff} mentioned 
in Section~\ref{s:dec}:
\begin{thm}
Let $\FF$ be a smooth codimension $q$ foliation on $M$ and $\Phi_t$ 
a smooth $1$-parameter family of diffeomorphisms of $M$. If $\FF_t$ is 
the family generated by pulling back $\FF$ via $\Phi_t$, 
then $T(\FF_t )=0\in H^{2q+1}(\II(\FF_{t}))$.
\end{thm}
\begin{proof}
In this situation the image of $T(\FF_t )$ in the de Rham cohomology is 
trivial because the flow acts trivially on de Rham cohomology, and so 
$T(\FF_t )$ is represented by an exact form. The point of the Theorem, 
and of this proof, is to see that the primitive can be chosen to be in 
the ideal $\II(\FF_{t})$.

Note that~\eqref{eq:dbe} shows that $(d\bb)^{q}=\al\ww\cc$, where 
$\cc$ is a locally decomposable $q$-form. Under the assumption of 
the Theorem, we may choose $\al$ so that $\dot\al =L_X\al$,
with $X$ a (time-dependent) vector field. 

Now we calculate:
\begin{alignat*}{1}
\dot\bb\ww (d\bb)^{q} \ &= \ \dot\bb\ww\al\ww\cc
\ = \ (-1)^{q}\al\ww\dot\bb\ww\cc\\
&= \ (-1)^{q}(d\dot\al)\ww\cc+ (-1)^{q-1}\dot\al\ww\bb\ww\cc\\
&= \ (-1)^{q}(di_{X}d\al)\ww\cc+ (-1)^{q-1}(di_{X}\al)\ww\bb\ww\cc\\
&\ \ \ \ +(-1)^{q-1}(i_{X}d\al)\ww\bb\ww\cc\\
&= \ (-1)^{q}(di_{X}(\al\ww\bb))\ww\cc+ 
(-1)^{q-1}(di_{X}\al)\ww\bb\ww\cc\\
&\ \ \ \ +(-1)^{q-1}(i_{X}(\al\ww\bb))\ww\bb\ww\cc\\
&= \ (-1)^{q}(d((i_{X}\al)\ww\bb+ (-1)^{q}\bb(X)\al))\ww\cc\\
&\ \ \ \ +(-1)^{q-1}(di_{X}\al)\ww\bb\ww\cc -\bb(X)\al\ww\bb\ww\cc\\
&= \ (d\bb(X))(d\bb)^{q}= d(\bb(X)(d\bb)^{q}) \ ,
\end{alignat*}
which is clearly in $d(\II(\FF_{t}))$. Here we have used 
$d\bb\ww\cc=0$, which follows from~\eqref{eq:dbe} and the definition 
of $\cc$ as the product (up to a constant) of the $\cc_{i}$.
\end{proof}
    
If $M$ is closed and oriented,
we can think of $TGV(\FF_{t})$ as a linear functional on 
$H^{n-2q-2}(M,\R)$. Dualising the above decomposition of 
$TGV(\FF_{t})$, we see that this functional factors through the 
Godbillon operator defined by multiplication with the Reeb 
class as in Theorem~\ref{t:dual}. The vanishing theorems for the 
Godbillon operator have been extended to higher codimension by 
Hurder~\cite{Hu}. Combining his result with Theorem~\ref{t:dvhigher} 
we obtain:
\begin{thm}
Let $\FF_{0}$ be a foliation almost all of whose leaves have 
subexponential growth. Then $DGV(\FF_{t})$ vanishes for all germs $\FF_{t}$ 
at $\FF_{0}$.
\end{thm}

\section{Final comments}\label{s:last}

In this paper we have extended the Godbillon-Vey invariants to 
families of foliations, obtaining families of cohomology classes in 
even degrees. This can be done quite generally, for all the 
characteristic classes of foliations, and is the subject of joint work 
with M.~Hoster and F.~Kamber~\cite{HKK}. 
The general construction starts from the observation that a 
$1$-parameter family $\FF_{t}$ of codimension $q$ foliations on $M$ 
can be thought of as a foliation of codimension $q+1$ on $M\times\R$, 
such that each cross-section $M\times\{t\}$ is saturated. 
At the level of forms, or in the foliated cohomology, $TGV(\FF_{t})$ 
is then obtained from the classical Godbillon-Vey invariant in 
codimension $q+1$ by contraction with the vector field 
$\frac{\partial}{\partial t}$. Of course, this is not a cohomological 
calculation in de Rham cohomology.
In the general case, however, the explicit calculations with differential 
forms become very complicated, and we resort to the 
formalism of Weil algebras. 
%This also allows us to make contact with 
%the work of Gelfand, Feigin and Fuks~\cite{GFF,Fuks}. 
Further generalisations to the case of multi-parameter families and 
to flags of foliations are discussed in Hoster's thesis~\cite{H}.

Heitsch~\cite{He} considered time-derivatives of characteristic classes, 
which are classes of the same degree as the classical characteristic 
classes of foliations. He showed that his classes vanish for all 
infinitesimal variations of the Roussarie example~\cite{GV}. As this 
is transversely homographic, Heitsch's result is a consequence of 
Theorem~\ref{t:homo} above. 

\medskip

{\sl Acknowledgements:} I am grateful to R.~Bott, M.~Hoster and 
F.~Kamber for many useful discussions. Hoster in particular provided 
valuable help with some of the calculations in Section~\ref{s:higher}.

\bibliographystyle{amsplain}

\bigskip

\end{document}